\documentstyle[cd,12pt]{article}
\newtheorem{thm}{Theorem}

\newtheorem{lemma}{Lemma}

\newtheorem{corol}{Corollary}

\newcommand{\T}{Teichm\"{u}ller }

\newcommand{\C}{$\hat {\Bbb C}$}

\newcommand{\Bbb}{\bf}
\newcommand{\pf}{{\it Proof}.~}
\newcommand{\qed}{$\hfill{\Box}$}

\newcommand{\cdr}{\cd\rightarrow}

\newcommand{\cdd}{\cd\downarrow}

\title{The Bers-Greenberg Theorem
and the Maskit Embedding for Teichm\"uller spaces}
\author{Pablo Ar\'{e}s Gastesi\thanks{\noindent Partially supported by SUNY at
Stony Brook and Tata Institute.\newline
}}

\begin{document}
\date{\today}
\maketitle

\begin{center}
{\it To the memory of Lipman Bers}
\end{center}

\begin{abstract}
The Bers-Greenberg theorem tells that the \T space of a Riemann
surface with branch points (orbifold)
depends only on the genus and the number of special points, but not on
the particular ramification values. On the other hand, the
Maskit embedding provides a mapping from the \T space of an
orbifold, into the product of one
dimensional \T spaces. In this paper we prove that there is a set of
isomorphisms between one dimensional \T spaces that, when restricted
to the image of the \T space of an orbifold under the Maskit
embedding, provides the Bers-Greenberg isomorphism.
\end{abstract}

\section{Introduction.}
The purpose of this paper is to connect two important results
in Teichm\"{u}ller theory, namely the Bers-Greenberg theorem
and the Maskit embedding.   We will prove that the isomorphism
given by the former theorem splits into the product of one
dimensional isomorphisms in the Maskit embedding.

Let $S$ be an orbifold with hyperbolic signature $(p,n;\nu_1,
\ldots, \nu_n)$.  A key result of B. Maskit tells that $S$ can be
written as the quotient $\Delta/\Gamma$, where $\Gamma$ is a
terminal regular b-group (it uniformizes $S$ and a finite number
of orbifolds of type $(0,3)$, which carry no moduli) and
$\Delta$ is a simply connected $\Gamma-$invariant open
subset of the regular set $\Omega(\Gamma)$ of $\Gamma$.  The
Bers-Greenberg theorem says that the complex structure of
the \T space of $\Gamma$, $T(\Gamma)$, depends only on
$(p,n)$, and not on the particular values $\nu_j$.

Applying the combination theorems (these are just tools to construct
Kleinian groups from smaller groups), Maskit produced an embedding
of the  Teichm\"{u}ller space $T(\Gamma)$ into the product
${\displaystyle{\prod_{j=1}^{3p-3+n}}} T(\Gamma^j),$ where the
$\Gamma^j$'s are certain subgroups of $\Gamma$, whose \T
spaces have dimension $1$. One can ask whether it is
possible to put the Bers-Greenberg theorem and the
Maskit embedding together in a single result.  More precisely, we
have a diagram
$$\CD
T(\Gamma_0) \cdr{h^*}{} T(\Gamma) \\
\cdd{i_0}{} \cdd{}{i} \\
\displaystyle\prod_{j=1}^{3p-3+n}T(\Gamma_0^j) \cd.
\displaystyle\prod_{j=1}^{3p-3+n}T(\Gamma^j),
\endCD$$
where $\Gamma_0$ and $\Gamma$ are terminal regular b-groups of
the same type $(p,n), i_0$ and $i$ are the Maskit embeddings,
and $h^*$ is an isomorphism obtained in the proof of the
Bers-Greenberg theorem.  The main result of this paper is that
the above diagram can be closed in a commutative way with a set
of one dimensional isomorphisms.
\begin{thm} Let $S$ be an orbifold with hyperbolic
signature.  Let $S_0$ be the orbifold obtained by removing from
$S$ all points with finite ramification value.  Let $\cal P$ be a
maximal partition on $S$ (and $S_0)$.  Assume $\Gamma$ and
$\Gamma_0$ are terminal regular b-groups uniformizing $(S,\cal P)$
and $(S_0,\cal P)$ on the invariant components $\Delta$ and
$\Delta_0$ respectively.  Then there exist a choice of
modular subgroups $\Gamma^1, \ldots, \Gamma^{3p-3+n}$ of
$\Gamma$, and $\Gamma_{0}^{1}, \ldots, \Gamma_{0}^{3p-3+n}$ of
$\Gamma_0$, and a set of isomorphisms
$h_{j}^{*}: T(\Gamma_{0}^{j}) \rightarrow T (\Gamma^j), 1 \leq j
\leq 3p-3+n,$ such that the restriction of $(h^{*}_{1}, \ldots,
h^{*}_{3p-3+n})$ to $i_0(T(\Gamma_0))$ closes the above diagram
in a commutative way.
\end{thm}
{\bf Remark}. Observe that we can choose $i^{-1} \circ h^{*}
\circ i_0$ as an isomorphism between $i_0 (T(\Gamma_0))$ and
$i(T(\Gamma))$, but our result is deeper
since it produces isomorphisms at the one dimensional level that
are ``compatible'' with the Maskit embedding.

This paper is organized as follows.  In $\S 2$ we give the
necessary background on Kleinian groups and Teichm\"{u}ller
spaces; in $\S 3$ we prove the Bers-Greenberg theorem adjusted
to our case following \cite{ek:hol}; $\S 4$ contains some
technical results; finally in $\S 5$ we prove theorem 1.

\noindent {\bf Acknowledgements.}  I would like to thank my advisor, Irwin
Kra, for posing the question that made possible this paper and
for all his help.  Special thanks are due to M. Lyvbich, who
provided the proof of lemma 4.3., and C.Zhang for many useful
conversations.  I am also very grateful to the referee for many
comments that helped to clarify this paper. Thanks are also due to K. P.
Sivaraman the typing the manuscript.

\section{Background on Kleinian groups and \T spaces.}
\noindent
{\bf 2.1.}  Let $\overline{S}$ be compact Riemann surface of
genus $p$, and let $x_1,\ldots, x_n$ be $n \geq 0$ distinct
points of $\overline{S}$.  Let $\nu_1, \cdots, \nu_n \in Z^+
\cup \{ \infty\}, \nu_j \geq 2,$ satisfying $2p-2+n -
\sum_{j=1}^{n} \frac{1}{\nu_j} >0$.  The tuple
$\sigma = (p,n;\nu_1,\ldots,\nu_n)$ is called a {\it
hyperbolic signature}, and $(p,n)$ is called the {\it
type} of $\sigma$; the $\nu_j$'s are called {\it ramification values}.
Assign the value $\nu_j$ to $x_j$.  An {\bf orbifold} structure on
$\overline{S}$ is given by a branched covering from the upper
half plane $H$ onto $S=\overline{S}-\{x_j; v_j= \infty\}$, such
that the natural projection $H \stackrel{\pi}{\rightarrow} S$ is
$\nu_j$-to-$1$ in a neighbourhood of $x_j$ whenever $\nu_j <
\infty$ (by an abuse of notation we will say that $S$ is an
orbifold with signature $\sigma$ and of type $(p,n)$).  It is a
consequence of the Poincar\'{e}-Koebe uniformization theorem that there
exists a Fuchsian group $F$, such that $H/F \cong S$, but we
are interested on a different representation of $S$.

A {\bf maximal partition}, ${\cal P}=\{a_1,\ldots, a_{3p-3+n}\}$, on $S$
is a set of simple closed disjoint unoriented curves on
$S_0=S-\{x_j;\nu_j<\infty\}$, such that no curve of $\cal P$ bounds
a disc or a punctured disc on  $S_o$, and no two curves of $\cal P$
bound a cylinder on $S_0$.

For a Kleinian group $G,$ let $\Omega (G)$ denote the non-empty
open set of points of $\stackrel{\wedge}{\bf C}=C\cup\{\infty\}$
where $G$ acts discontinuosuly. The components of $\Omega (G)$ are called
the components of $G$.
\begin{thm} [Maskit uniformization] Given an orbifold $S$
with hyperbolic signature $\sigma$ and maximal partition $\cal P$, there
exists a (unique up to conjugation  by a M\"{o}bius
transformation) geometrically finite Kleinian group $\Gamma$,
called a {\bf terminal regular b-group}, with a simply connected
invariant component $\Delta$, such that:

1.- $\Delta / \Gamma$ is conformally equivalent to $S$;

2.-  the curves of $\cal P$ are in one-to-one correspondence with the
conjugacy classes of maximal cyclic subgroups of $\Gamma$
generated by accidental parabolic elements;

3.- $(\Omega (\Gamma) - \Delta)/\Gamma$ is the (finite) union
of the orbifolds of type $(0,3)$ obtained by squeezing each
curve of $\cal P$ to a point, and discarding all orbifolds of
signature $(0,3,2,2,\infty)$ that appear.
\end{thm}
See Maskit's book \cite{mas:kg} for the definition of geometrically finite and
accidental parabolic elements; they are technical terms not needed to
understand the proofs of this paper.

\smallskip\noindent
{\bf 2.2.} Let $G$ be a finitely generated non-elementary
Kleinian group, and let $A$ be an open $G-$invariant subset of
$\Omega (G)$.  Let $L^{\infty}(G,A)$ denote the set of
essentially bounded (classes of) functions $\mu$, supported on
$A$ and such that $(\mu\circ\gamma) \overline{\gamma '} / \gamma '
=\mu$, for all $\gamma \in G$.  Let $M(G,A)$ be the unit
ball of $L^{\infty} (G,A);$ the classical theory of
quasiconformal mappings tells us
that if $w$ is a quasiconformal homeomorphism of
$\stackrel{\wedge}{{\bf C}}$ satisfying the  equation $(*)~
w_{\overline{z}} = \mu w_z$, for some $\mu \in M(G,A)$,
then the group $w \Gamma w^{-1}$ is again a group of M\"{o}bius
transformations (see \cite{nag:teic}).  Two functions in $M (G,A)$,
$\mu_1$ and $\mu_2$,
are {\it equivalent} if $w^{\mu_1} \gamma (w^{\mu_1})^{-1} =w^{\mu_2}
\gamma (w^{\mu_2})^{-1}$ for all $\gamma \in G,$ where
$w^{\mu}$ is the unique solution of $(*)$ fixing $\infty$, $0$ and $1$
(as a technical assumption, we need to suppose that these three
points do not lie in $\Omega (G);$ this can always be achieved
by conjugating by a M\"{o}bius transformation since $G$ is
non-elementary).  The set of equivalence classes of elements of
$M(G,A)$ is called the {\bf Teichm\"{u}ller} space of $G$
supported on $A$, and it is denoted by $T(G,A)$.  This set has a
unique complex structure so that the natural projection from
$M(G,A)$ onto $T(G,A)$ is holomorphic.  In the case we are
interested on we will take $G=\Gamma$ as in the Maskit
uniformization theorem, and $A=\Delta$; we then have
$T(\Gamma,\Delta)= T(\Gamma,\Omega (\Gamma))$ (since the
orbifolds uniformized by $\Gamma$, other than $S$, carry no
moduli). We will denote this set by $T(\Gamma )$.
It is a well known fact that $T(\Gamma)$ is a complex
manifold of dimension $3p-3+n$.  One can ask up to what point
the complex structure of this space depends on the values
$\nu_j$ appearing in the signature $\sigma$ of $S
\cong \Delta/\Gamma$.  The following theorem
answers that question.
\begin{thm} [Bers-Greenberg \cite{bg:isom}] Assume $\Gamma_1$ and
$\Gamma_2$ are two terminal regular b-groups with invariant
components $\Delta_1$ and $\Delta_2$ such that the
orbifolds $\Delta_1/ \Gamma_1$ and $\Delta_2/ \Gamma_2$
are of the same type $(p,n)$.  Then the Teichm\"{u}ller spaces
$T(\Gamma_1)$ and $T(\Gamma_2)$ are biholomorphically
equivalent.
\end{thm}
{\bf Remark}.  Although the theorem was initially proven for
Fuchsian groups, the result is also valid in the case of
terminal regular b-groups; we will provide a proof of it in
$\S 3$.

\bigskip\noindent
{\bf 2.3.}  The Maskit uniformization theorem
produces an embedding of $T(\Gamma)$ into the product
of one-dimensional Teichm\"{u}ller spaces as follows.  Let $T_j$,
$1\leq j\leq 3p-3+n$,
be the connected component of $S-\{a_k; a_k \in C, k \neq
j\}$ containing $a_j$ (these sets are called the {\it modular}
parts of $S)$.  Let $D_j$ be a component of $\pi^{-1} (T_j)$, where
$\pi:\Delta \rightarrow S$ is the natural projection, and let
$\Gamma^j=\{\gamma \in \Gamma; \gamma (D_j)=D_j\}$.  The
groups $\Gamma^j$ are known as {\it modular} subgroups of
$\Gamma$; observe that if $\tilde{D}_j$ is another component of
$\pi^{-1} (T_j)$ and $\tilde{\Gamma}^j $ is the corresponding
subgroup, then $\Gamma^j$ and $\tilde{\Gamma}^j$ are conjugated
on $\Gamma$.  From the work of Maskit one can see that
$\Gamma^j$ is a terminal regular b-group of type $(0,4)$ or
$(1,1)$ with invariant component $\Delta^j \supset
\Delta$.  It is clear that $M (\Gamma^j, \Delta^j)$
contains $M(\Gamma, \Delta)$, but we have a stronger result.
\begin{thm} [Maskit embedding;
Maskit \cite{mas:moduli}, Kra \cite{kra:variational}]
The mapping given by restriction $T(\Gamma) \rightarrow
{\displaystyle{\prod_{j=1}^{3p-3+n}}} T(\Gamma^j)$ is
holomorphic and injective with open image.
\end{thm}

By a result of Gentilesco \cite{gent:aut}, the image of $T(\Gamma)$ is not
the whole set $\prod_{j=1}^{3p-3+n} T(\Gamma^j)$,
unless we are in the trivial case of dim
$T(\Gamma)=1$.

\bigskip\noindent
{\bf 2.4.}  We need to define one more set of holomorphic
functions related to Kleinian groups.  Given $\Gamma$ and
$\Delta$ as above (theorem 2), $Q(\Gamma, \Delta)$ will
denote the space of {\it quadratic differentials} (for $\Gamma$
on $\Delta)$ consisting on the functions $f$, holomorphic on
$\Delta$ such that $(f\circ \gamma)(\gamma ')^{2}=f$ for all $\gamma
\in \Gamma$, and with finite norm
$$|| f || = \frac{1}{2} \int
\int_{\Delta / \Gamma} \mid f(z) dz \wedge d \overline{z} \mid
< \infty.$$
The elements of $M(\Gamma, \Delta)$ of the form $k
\frac{\overline{\varphi}}{\mid \varphi \mid}$, for some $\varphi
\in Q (\Gamma, \Delta)$ and some $k \in (0,1)$
real, are called {\bf Teichm\"{u}ller differentials}.
Teichm\"{u}ller's  theorem (Ahlfors \cite{ahl:qc}, Bers \cite{bers:teichm})
tells that on each
class of $T(\Gamma)$ there exists a unique Teichm\"{u}ller
differential.

\bigskip\noindent
{\bf 2.5.} Finally we define the Teichm\"{u}ller space of an
orbifold and study the relationship with the Teichm\"uller space
of Kleinian groups.  Let $S$ be an orbifold with hyperbolic
signature $\sigma$; consider the  set of quasiconformal homeomorphisms
({\it deformations}) $f:S \rightarrow S'$, where $S'$ is
another orbifold with the same signature than $S$, and such that
the ramification values of $x$ and $f(x)$ are the same for all
$x \in S$.  Two such mappings $f:S \rightarrow S'$ and
$g:S \rightarrow S"$, are {\it equivalent} if there exists a
biholomorphic function $\phi:S' \rightarrow S"$ (respecting the
ramification values) such that $g^{-1} \circ \phi \circ f$ is
homotopic to the identity on $S$ (by a homotopy that fixes the
points $x_j)$.  The set of equivalence classes of deformations
of $S$ is the {\bf Teichm\"{u}ller space} of $S$ ,$T(S)$.  We have
that if $S \cong \Delta / \Gamma$, with $\Gamma, \Delta$
as in theorem 2, or $S \cong H/F$ with $F$ Fuchsian, then $T(S)
\cong T(\Gamma)$ and $T(S) \cong T(F,H)$ (usually denoted by
$T(F))$; see \cite{kra:spaces} and \cite{nag:teic}.

\section{The Bers-Greenberg isomorphism for b-groups.}
In \cite{ek:hol} I. Kra gave a proof of the Bers-Greenberg theorem for
the case of Fuchsian groups.  In this section, we will follow
his arguments to provide a proof of the theorem in the
case of terminal regular b-groups.  See also the remark at the
end of the proof.

\bigskip\noindent
{\it Proof of theorem 3.} Let $S$ be an orbifold with hyperbolic
signature $\sigma$, and assume that at least one of the ramification
values of $\sigma$ is finite.  Let $\cal P$ be a maximal partition on
$S$.  By theorem 2 we have that the pair $(S,\cal P)$ can be
uniformized by a terminal regular b-group $\Gamma$ acting on its invariant
component $\Delta$.  Remove from $S$ the points with finite
ramification value to obtain a surface with punctures $S_0$;
applying theorem $2$ again, we get a terminal regular b-group
$\Gamma_0$ uniformizing $(S_0,\cal P)$ on the invariant component
$\Delta_0$.  Consider now the set $\Delta_{\Gamma} =
\Delta - \{$ fixed points of elliptic elements of $\Gamma\}$.
Since we are assuming that the signature of $S$ contains finite
ramification values, we have $\Delta_{\Gamma} \neq \Delta$.
We also have $\Delta_\Gamma/ \Gamma \cong S_0$, so there exists a
holomorphic covering $h: \Delta_0 \rightarrow \Delta_r$
making the following diagram commutative:
$$\CD
\Delta_0 \cdr{h}{} \Delta \\
\cdd{\pi_0}{} \cdd{}{\pi} \\
S_0 \cdr{id}{} S_0
\endCD$$
where $\pi_0$ and $\pi$ are the natural projections.

$h$ induces a group homomorphism $\chi:\Gamma_0 \rightarrow
\Gamma$ as follows. For any  $\gamma \in \Gamma_0$, the
functions $h$ and $h \circ \gamma$ are lifts of $id \circ \pi_0$
by $\pi$, and therefore there exists an element $g \in
\Gamma$ such that $h \circ \gamma = g \circ h;$ define $\chi
(\gamma)= g$.

Using the map $h$ we can define a norm-preserving isomorphism
$h^*:L^{\infty} (\Gamma_0, \Delta_0) \rightarrow L^{\infty}
(\Gamma, \Delta)$, given by $(h^* \mu) \circ h=\mu h'/
\overline{h'}$.  One can check (see lemma after the proof) that
$h^*$ induces a holomorphic mapping between $T(\Gamma_0)$
and $T(\Gamma)$ (which
we will also denote by $h^*$).

$h$ also gives a mapping between quadratic differentials,
$h_*:Q (\Gamma, \Delta) \rightarrow Q(\Gamma_0, \Delta_0)$
by $h_* \varphi = (\varphi_0 h)(h')^2$.  The formula
$$h^* (\frac{\overline{h_* \varphi}}{\mid h_* \varphi \mid})=
\frac{\overline{\varphi}}{\mid \varphi \mid},$$
for $\varphi \epsilon Q(\Gamma, \Delta),$ shows that $h^*$
takes Teichm\"uller differentials to Teichm\"uller differentials
and therefore $h^*$ is bijective by Teichm\"uller's theorem.\qed

\bigskip\noindent
{\bf Remark}.  One can provide a shorter proof of the above
theorem as follows.  Given $S$ and $S_0$ as above, find Fuchsian
groups $F$ and $F_0$ such that $H/F \cong S$ and $H/F_0 \cong S_0$.
We then have $T(F) \cong T(S) \cong T(\Gamma)$ and $T(F_0) \cong
T(S_0) \cong T(\Gamma_0)$. By \cite{bg:isom} and \cite{ek:hol},
we know that $T(F)\cong T(F_0)$, and therefore $T(\Gamma)\cong T(\Gamma_0)$.
But for the purpose of this paper we
need to have a concrete expression of such an isomorphism, namely
$h^*$.
\begin{lemma} $h^*:T(\Gamma_0) \rightarrow T(\Gamma)$ is well
defined.
\end{lemma}
\pf It suffices to consider the case of $\mu$
being equivalent to $0$, i.e. $w \gamma w^{-1}=\gamma$, for all
$\gamma \in \Gamma_0$, where $w=w^{\mu}$.  In that case we
get $w (\Delta_0) =\Delta_0$ since this set is the unique
invariant component of $\Gamma_0$.  Define a mapping
$f:\Delta_\Gamma \rightarrow \Delta_\Gamma$ by $f(h(z))=h(w(z))$. By
considering the exact sequence
$$\{1\} \rightarrow ker \chi \rightarrow \Gamma_0
\stackrel{\chi}{\rightarrow} \Gamma \rightarrow \{1\},$$
it is not hard to see that $f$ is well defined and one-to-one.
So we get that $f$ is a quasiconformal mapping whose
coefficient is $h^*(\mu)$.  Extending this coefficient to \C ~by $0$
outside $\Delta_\Gamma$, we
get $f=A \tilde{w}$, where $A$ is a M\"obius transformation and
$\tilde{w} = w^{h^{*}(\mu)}$.

\noindent Let $g \in \Gamma$ and let $\gamma \in \Gamma_0$ be
such that $\chi (\gamma)=g$.  We then have

$$fgf^{-1} h=fghw^{-1}= fh \gamma w^{-1}=hw \gamma w^{-1} =h\gamma,$$
which implies that $\chi (\gamma)=fgf^{-1}$, that is
$fgf^{-1}=g$.  In particular we have that $A \tilde{w} (x) =x$,
for all $x \in \partial \Delta$; since we are assuming
that $\infty$, $0$ and $1$ lie in $\partial \Delta$, we get $A=id$
and therefore $\tilde{w}$ must be trivial.\qed

\section{Some technical results.}
In this section we will prove some technical lemmas needed in the
proof of theorem 1.  We deal mainly with three points: first
we show that the homomorphism $\chi$ of $\S 2$ takes modular subgroups
of $\Gamma_0$ onto modular subgroups of $\Gamma$.
The second point is a lemma about
Teichm\"{u}ller spaces (deformation lemma) that has interest on
its own, and finally we prove some properties of the Maskit
embedding.

\bigskip\noindent
{\bf 4.1.} Let $S,~S_0,~\Gamma,~\Gamma_0,~\Delta$ and
$\Delta_0$ be as in $\S 2$.  Let $T$ be one of the modular
parts of $S_0$, and let $\pi^{-1} (T)={\displaystyle{\cup_{j \in
J}}} D_j$ be a decomposition of the pre-image of $T$ into
disjoint connected components (here $\pi$ is the natural
projection of diagram $2$).  Apply $h^{-1}$ and
the commutativity of the diagram $2$ to get
$$\pi_0^{-1}(T) = h^{-1}(\pi^{-1} (T))=
{\displaystyle{\cup_{j \in J}}} h^{-1} (D_j).$$
On the other hand, we also have a decomposition into connected
components given by $\pi_0^{-1}(T) = {\displaystyle{\cup_{k
\in K}}} A_k$.  Therefore, for all $k \in K$ there
exists a $j \in J$ such that $h(A_k) \subset D_j$.
\begin{lemma} $h(A_k)= D_j$.
\end{lemma}
\pf Let us see first that if $A_k$ and $A_l$ are
such that $h(A_k)\subset D_j$ and $h(A_l)\subset D_j$, then $h(A_k)
\cap h (A_l)=\emptyset$ or $h(A_k)=h(A_l)$.  So assume $y$
belongs to $h (A_k) \cap h (A_l),$ and let $x_1\in
A_k,~ x_2 \in A_l$ be such that $h(x_1)=y=h(x_2)$.
This gives $\pi_0 (x_1)= \pi_0 (x_2)$, so there is a
$\gamma \in \Gamma_0$ such that $\gamma (x_1)=x_2$.  In
particular we get $\gamma (A_k) \cap A_l \neq
\emptyset$.  Since $\gamma$ is a homeomorphism and $\pi_0 \circ
\gamma = \pi_0$, we have ${\displaystyle{\cup_{p \in K}}}
A_p= {\displaystyle{\cup_{p \in K}}} \gamma (A_p)$, so
$\gamma (A_k)=A_l$.  We also obtain the equalities
$$y=h(x_2)= h (\gamma (x_1))=\chi (\gamma) (h (x_1))= \chi
(\gamma) (y).$$
Since the only transformation of
$\Gamma$ with fixed points in $\Delta_\Gamma$
is the identity, we must have
$\chi (\gamma)= id$, and therefore
$$h(A_l)= h (\gamma (A_k))= \chi (\gamma) (h (A_k))= h
(A_k)$$
as claimed.

Now to complete the proof of the lemma, let $L$ be
the set $L=\{k \in K; h (A_k) \subset D_j\}$.  Pick any
$k_0 \in L$ and divide $L$ into two disjoint sets, $L=M
\cup N$, as follows:  $M=\{k \in L; h (A_k)=h (A_{k_0})\}$ and
$N=\{k \in L; h (A_k) \cap h (A_{k_0}) = \emptyset \}$.  Then $D_j=
{\displaystyle{\cup_{k \in L}}} h(A_k)=
{\displaystyle{\cup_{k \in M}}}  h (A_k) \cup W$, where
$W={\displaystyle{\cup_{k \in N}}} h (A_k)$.  The sets
$h(A_{k_0})$ and $W$ are open and disjoint, so $W=\emptyset$ (since
$D_j$ is connected), that is, $N=\emptyset$, giving $h(A_{k_0})=D_j$.
\qed

\bigskip\noindent
{\bf 4.2.} Choose $D_j$ and $A_k$ with $h(A_k)=D_j$ and rename
$A_j=A_k$ for simplicity.  Let $\Gamma_0^{j}=\{\gamma \in
\Gamma_0; \gamma (A_j)=A_j\}, \Gamma^j=\{\gamma \in \Gamma; \gamma
(D_j)=D_j\}$.  It is easy to see that $\chi (\Gamma_0^{j}) \subset
\Gamma^j$, where $\chi$ is the group homomorphism of $\S 2$, but
we have a stronger result.
\begin{lemma} $\chi (\Gamma_0^{j})=\Gamma^j$.
\end{lemma}
Consider the diagram:
$$\CD
A_j \cdr{h}{} D_j \\
\cdd{\pi_0}{} \cdd{}{\pi} \\
T_0 \cdr{id}{} T_0
\endCD$$
Let $g \in \Gamma^j$ and let $z \in D_j$.  Since $h
(A_j)=D_j$, there exist points $x,x' \in A_j$ such that
$h(x)=z$ and $h(x')=g(\gamma)$.  We then get $\pi_0 (x')=\pi
(h(x'))=\pi (g(z))= \pi (z) = \pi (h(x))= \pi_0 (x')$, so there is
$\gamma \in \Gamma^j$ such that $\gamma (x) =x'$.   Now
the functions $h\circ \gamma$ and $g\circ h$ satisfy
$\pi \circ h\circ  \gamma=\pi\circ g\circ h$.
So to prove the equality $h\circ \gamma=g\circ h$ all we need to do is
to check it at one point.  Choose $x$ as before; then $h (\gamma (x))=h(x')=
g(z)$ and $g(h (x))=g(z)$, which proves that $\chi (\gamma )=g$ and therefore
the lemma.\qed

\bigskip\noindent
{\bf 4.3.} Consider now one of the groups we are working with,
say $\Gamma$, and a modular subgroup of it, $\Gamma^j$.  Since
$\Delta$ is invariant under $\Gamma$, it will also be invariant under
$\Gamma^j$, and therefore we can look at the Teichm\"uller space
of $\Gamma^j$ supported on $\Delta, T(\Gamma^j,
\Delta)=\{[\mu] \in T^j;$ supp $(\mu) \subset
\Delta\}$, where $[\mu]$ denotes the equivalence class of $\mu$.  The
following lemma says that this set is actually the whole
$T(\Gamma^j)$.
\begin{lemma} [Deformation lemma. Analytic version]
$T(\Gamma^j, \Delta)=T(\Gamma)$.
\end{lemma}

\noindent Before giving the proof of
this lemma, let us look at its geometric meaning.  The group
$\Gamma^j$ is a terminal regular b-group of type $(1,1)$ or
$(0,4)$, with invariant component $\Delta^j$.  The proof is
essentially the same in all cases, so we will assume that
$\Gamma^j$ has signature $(1,1; \infty),$ that is, $\Delta^j/
\Gamma^j$ is a torus with one puncture.  We have that the
quotient $\Delta/\Gamma^j$ is a torus with a hole. The
difference $\Delta^j/\Gamma^j - \Delta/
\Gamma^j$ is a punctured disc $D$.  If $\mu \in M (\Gamma^j,
\Delta^j)$ is supported on $\Delta$, the corresponding
homeomorphism induced on $\Delta^j/ \Gamma^j$ will be
conformal on $D$ (recall that $\mu=0$ implies $w^{\mu}$ is
conformal).  Therefore the above lemma is equivalent to the
following statement.
\begin{lemma} [Deformation lemma. Geometric version.]  Let $T$
be a surface with signature $(1,1; \infty)$, and let $D$ be a
punctured disc on $T$ containing the puncture.  Then any
quasiconformal deformation of $T$ is equivalent to a
deformation that is conformal on $D$.
\end{lemma}
\pf Let $f:T\rightarrow T'$ be a quasiconformal homeomorphism.
Let $D'=f(D)$; since $f$ is quasiconformal, $D'$ is a also punctured disc.
Any two punctured discs are conformally
equivalent, so there exists a conformal homeomorphism $g:D
\rightarrow D'$.  Let $h$ denote the function $g \circ
f^{-1}:D' \rightarrow D'$.  Let $V'$ be a punctured disc on
$T'$ containing $D'$ and such that $V'- D'$ is annulus.
Since $h$ is an orientation preserving homeomorphism, we can
extend it smoothly to $T'$ so that (the extension of) $h$ is the
identity on $T'-V'$.  Now the functions $f$ and $h
\circ f$ are homotopic, therefore equivalent, and $h \circ f=g$
is conformal on $D$.\qed

\bigskip The above proof shows that the result can be generalized to any
hyperbolic orbifold as follows.
\begin{lemma} [Deformation Lemma.] Suppose $S$ is an orbifold
with hyperbolic signature $(p,n; \nu_1, \ldots, \nu_n)$. Let
$D_1, \ldots, D_n$ be discs around the special points of $S$ (or
punctured discs around the punctures) such that their closures
are pairwise disjoint.  Then any quasiconformal deformation of
$S$ is equivalent to a deformation which is conformal on $U_1\cup
\ldots \cup U_n$.
\end{lemma}

\bigskip\noindent
{\bf 4.4.} Our last lemma gives some properties of the Maskit Embedding for
\T spaces of b-groups. As in lemma 4.3., all the proofs are independent
(up to technical points) of the ramification values, so for the
sake of simplicity we will work out the torsion free situation,
that is, the case of b-groups uniformizing compact surfaces with
(possibly $0$) punctures.

Let us start by looking at the one-dimensional cases, namely those
of four times punctured spheres and once punctured tori. A b-group
$\Gamma$ of signature $(0,4;\infty,\ldots,\infty)$ is constructed by
taking two triangle groups, $\Gamma_1$ and $\Gamma_2$, of signatures
$(0,3;\infty,\infty,\infty)$, such that $\Gamma_1\cap\Gamma_2=<A>$, where
$A$ is a parabolic element (conjugate in $PSL(2,\Bbb C)$ to a translation);
then $\Gamma$ is the group generated by $\Gamma_1$ and $\Gamma_2$,
$\Gamma=\Gamma_1*_{<A>}\Gamma_2:=<\Gamma_1,~\Gamma_2>$. This is the
so-called AFP construction, and the First Combination Theorem
\cite[Theorem VII.C.2, pg. 149]{mas:kg} guarantees that $\Gamma$ is
Kleinian if we choose $\Gamma_1$ and $\Gamma_2$ properly. To fix ideas,
let $\Gamma_1$ be the group generated by $A(z)=z+2$ and
$B(z)=(-z)/(2z-1)$. Then $\Gamma_2=T_\alpha\Gamma_1T_\alpha^{-1}$,
where $T_\alpha(z)=z+\alpha$, with $\alpha\in\Bbb C$. If $Im(\alpha)$
is big enough, then we can take the horocircle
$\{z;~Im(z)=\frac{1}{2}Im(\alpha)\}$ as the invariant curve required in
the above quoted theorem, and $\Gamma$ is a b-group of the desired
signature. See figure $1$. In \cite{kra:horoc} (see \cite{ares:horoc1}
for the case of torsion. See also
\cite{km:defor}) it is proven that $\alpha$ is a global coordinate in
the space $T(\Gamma)\cong T(0,4;\infty,\ldots,\infty)$. One can write
$\alpha$ as the cross ratio $\alpha =cr(\infty,0,1,\alpha)$. Observe
that the four points involved in this expression are the fixed points of
the transformations $A$, $B$, $AB$ and $T_\alpha BT_\alpha^{-1}$,
respectively. This way of expressing $\alpha$ as a cross ratio has
the advantage of being independent of the particular choice of
$\Gamma_1$ and $\Gamma_2$; i.e., if we start with $D\Gamma_1D^{-1}$, for
some M\"obius transformation $D$, we still have a global coordinate in \T
space expressed as a cross ratio. We can see from the above discussion that
the set $\{z;~Im(z)>k\}$ is contained in $T(\Gamma)$ if $k$ is big emough
(actually, $k=1$ works).

For the case of surfaces with signature $(1,1;\infty)$, we
start with $\Gamma_1$ as above, and then we look for a M\"obius transformation
$C$ such that $CB^{-1}C^{-1}=A$ (see \cite{kra:horoc} for an explanation
of why we must take $B^{-1}$ and not $B$). This gives $C(z)=\tau+\frac{1}{z}$,
$\tau\in\Bbb C$. The horocircle
$\{z;~|z-\frac{ir}{2}|=\frac{r}{2}\}$ is mapped by $C$ to the horocircle
$\{z;~Im(z)=Im(\tau)-\frac{1}{r}\}$. These two curves will be disjoint
if $Im(\tau)>l$, for some positive constant $l$.
Then we can apply the
Second Combination Theorem \cite[VII.E.5, pg. 161]{mas:kg}
and obtain that the group
$\Gamma'=\Gamma_1*_C:=<\Gamma_1,C>$ is a b-group with the desired
signature. $\tau$ is a global coordinate on $T(\Gamma')=T(1,1;\infty)$, and
it can be expressed as the cross ratio $\tau=(\infty,0,1,C(\infty))$.
This construction is called an HNN extension. See figure $2$.

The general b-group $G$ uniformizing a surface with signature
$(p,n;\infty,\ldots,\infty)$ is constructed from triangle groups
(actually, $\Gamma_1$, its conjugates in $PSL(2,\Bbb C)$ and transformations
satisfying relations like $C$ above) by iterated applications of the
Combination Theorems. Suppose that we have done the first step and
constructed a b-group $G_1$ (of type $(0,4)$ or $(1,1)$). The following
construction will be either an AFP or an HNN extension. In either
case, if we use horocircles \lq close\rq~ to the fixed points of the
parabolic elements involved in the process, we are guaranteed that the
Maskit Theorems can be applied. But by the cross ratio expression of the
coordinates, one sees that these horocirlces correspond to points in
\T spaces with big imaginary parts. So we have proven the following
lemma, modulo some technical points that are required for the case of torsion,
but they do not present any difficulty.
\begin{lemma}  For any $\alpha \in T(G_i)$ there exist
non-negative real numbers, $y_{\alpha}^{2}, \ldots ,
y_{\alpha}^{3p-3+n}$, such that the  set $\{ \alpha\} \times
V_{\alpha}$ is contained in $i(T(G))$, where $V_{\alpha}
=\{(z_2, \ldots, z_{3p-3+n} \in C^{3p-4+n}; Im)(z_j)>
y_{\alpha}^{j}, j=2, \ldots, 3p-3+n\}$.
\end{lemma}
{\bf Remark.} It is clear that this lemma works for any $\alpha_j$
in $T(G_j), j=1, \cdots, 3p-3+n$.
\begin{corol}  Given $\alpha$ and $\beta$ in $T(G_1)$, we have
$V_{\alpha}\cap V_{\beta} \neq \emptyset$, where $V_\alpha$ an $V_\beta$ are
given by the above lemma.
\end{corol}

\section{Proof of theorem 1.}
Let $S, S_0, C_1, \Gamma, \Gamma_0, \Delta$ and $\Delta_0$
be as in $\S \S$ 2,3.  Choose a modular part $T_j$ of $S_0, 1
\leq j \leq 3p-3+n$. By lemma $1$ we can choose components $D_j$ and
$A_j$ of $\pi^{-1} (T_j)$ and $\pi_{0}^{-1} (T_j)$ respectively,
such that $h(A_j)=D_j$.  By lemma $2$ we have that $\chi
(\Gamma_{0}^{j})= \Gamma^j$, where $\Gamma_{0}^{j}$ and
$\Gamma^j$ are defined in $4.2.$,
and $\chi:\Gamma_0 \rightarrow
\Gamma$ is the group homomorphism of $\S$2.  Therefore we can
define a function $h^{*}_{j}$ by $(h_{j}^{*} \mu) \circ h = \mu
h'_j/ \overline{h'_j}$, where $h_j$ is the restriction of $h$ to
$A_j$, and $\mu$ belongs to $T(\Gamma_{0}^{j}, \Delta_0)$.
By the deformation lemma we actually  have that $h^{*}_{j}$ is
defined on $T(\Gamma_{0}^{j})$.  So we obtain mappings,
$h_{j}^{*}:T(\Gamma_{0}^{j}) \rightarrow T(\Gamma^j), 1 \leq j
\leq 3p-3+n$.  Since these mappings are defined by the same
differential equation as $h^*$, it is clear that the following
diagram commutes, if $f$ is equal to the restriction of
$(h^{*}_{1}, \cdots, h^{*}_{3p-3+n})$ to $i_0 (T(\Gamma_0)):$
$$\CD
T(\Gamma_0) \cdr{h^*}{} T(\Gamma) \\
\cdd{i_0}{} \cdd{}{i} \\
\displaystyle\prod_{j=1}^{3p-3+n}T(\Gamma_0^j) \cdr{f}{}
\displaystyle\prod_{j=1}^{3p-3+n}T(\Gamma^j),
\endCD$$
To complete the proof of theorem 1, all we need to show is that
each $h^{*}_{j}$ is bijective.  Take $j=1$ to simplify notation.
\vskip 5mm
\noindent
{\bf Injectivity}: Let $\alpha, \beta$ be in $T(\Gamma_{0}^{I})$
with $h_{1}^{*} (\alpha)= h_{1}^{*} (\beta)$.  Lemma $7$ gives two
open sets, $V_{\alpha},V_{\beta} \subset C^{3p-4+n}$ such that
$\{\alpha\} \times V_{\alpha}$ and $\{\beta\} \times V_{\beta}$
are contained in $i_0 (T(\Gamma_0))$.  By the corollary $1$ of the
lemma $7$ we have $V_{\alpha} \cap V_{\beta} \neq \emptyset$, so we can
choose $\gamma$ in that intersection.  We then get $f(\alpha,
\gamma)= f(\beta, \gamma)$, contradicting the injectivity of
$h^*$.
\vskip 5mm
\noindent
{\bf Surjectivity}:  Take $\alpha' \in T (\Gamma^1)$, and
let $\beta'$ be such that  $(\alpha', \beta') \in i(T(\Gamma))$
(which is possible by lemma $6$).  Then there exists a
point $x= (\alpha, \beta)$ in $i(T(\Gamma_0)) \subset
T(\Gamma_{0}^{1})\times \prod_{j=2}^{3p-3+n}
T(\Gamma_{0}^{j})$ such that $h^* (\alpha,\beta)=(\alpha',
\beta')$,  giving $h^{*}_{1} (\alpha)=\alpha'$. \qed

\ifx\undefined\bysame
\newcommand{\bysame}{\leavevmode\hbox to3em{\hrulefill}\,}
\fi

\newpage
\begin{picture}(250,200)(-150,0)
\thicklines
\put(-100,0){\line(50,0){250}}
\put(-100,200){\line(50,0){250}}
\put(-25,25){\line(0,50){150}}
\put(75,25){\line(0,50){150}}
\put(-25,0){\oval(50,50)[tr]}
\put(25,0){\oval(50,50)[tl]}
\put(25,0){\oval(50,50)[tr]}
\put(75,0){\oval(50,50)[tl]}
\put(-25,200){\oval(50,50)[br]}
\put(25,200){\oval(50,50)[bl]}
\put(25,200){\oval(50,50)[br]}
\put(75,200){\oval(50,50)[bl]}
\thinlines
\put(-100,100){\line(50,0){250}}
\end{picture}
\newline
\begin{center}Fig $1$. The AFP construction. \end{center}\vspace{7mm}
\begin{picture}(250,100)(-150,0)
\thicklines
\put(-100,0){\line(50,0){250}}
\put(-25,25){\line(0,50){95}}
\put(75,25){\line(0,50){95}}
\put(-25,0){\oval(50,50)[tr]}
\put(25,0){\oval(50,50)[tl]}
\put(25,0){\oval(50,50)[tr]}
\put(75,0){\oval(50,50)[tl]}
\thinlines
\put(-100,100){\line(50,0){250}}
\put(0,12.5){\circle{25}}
\end{picture}
\newline
\begin{center}Fig. $2$. The HNN extension construction. \end{center}

\noindent Math Department, SUNY at Stony Brook, USA and School of Maths,
Tata Institute of Fundamental Research, Bombay, INDIA.\\
\noindent {\it E-mail address}: pablo@motive.math.tifr.res.in

\end{document}